\documentclass{article}

\usepackage[all]{xy}
\usepackage[pdftex]{graphicx}

\usepackage[all]{xy}

\usepackage{amsmath,amssymb}
\usepackage{amsthm}

\newtheorem{prop}{Proposition}[section]

\theoremstyle{definition}
\newtheorem{defini}[prop]{Definition}

\newcommand{\ep}{\varepsilon}

\newcommand{\RR}{\mathbb R}

\newcommand{\ZZ}{\mathbb Z}
\newcommand{\QQ}{\mathbb Q}

\newcommand{\Crit}{{\rm Crit}}
\newcommand{\grad}{{\rm grad}}

\newcommand{\ind}{{\rm ind}}

\usepackage{url}

\begin{document}


\vspace*{1.5pc}

\large

{\bf\LARGE \begin{center}
A Diagrammatic method for an invariant related to the Chern-Simons perturbation theory
\end{center}}

\medskip

\centerline{\Large 
Tatsuro Shimizu
\footnote{School of System Design and Technology,
Tokyo Denki University\\
{\tt\hskip5mm t\underline{ }shimizu@mail.dendai.ac.jp}
}}

\bigskip
\begin{abstract}
The defect $d(M,\rho)$ is an invariant of a compact oriented 3-manifold $M$ with a representation $\rho$ of the fundamental group. 
In this article we give a diagrammatic method for $d$ of knot exteriors by using knot diagrams.
\end{abstract}

\medskip
\section{Introduction}
The defect $d$ is an invariant of a compact 3-manifold with a representation of the fundamental group. 
The defect $d$ can be found out in the construction of the Chern-Simons perturbation theory established by M.~Kontsevich \cite{Kon}, S.~Axelrod and I.~M.~Singer \cite{AS}. The defect $d$ was first formulated and studied by C. Lescop \cite{Lescop} for a closed oriented manifold $M$ with the first Betti number one and an abelian representation $\rho$. (In \cite{Lescop}, $d(M,\rho)$ is denoted by $I_{\Delta}(M)$). 
In the Chern-Simons perturbation theory, $d$ plays a role as a kind of defect.
Recently, H.~Kodani and B.~Liu in \cite{KodaniLiu} clarify a partial role of $d$ for adjoint representations of semi-simple Lie groups.

It is expected that $d(M,\rho)$ essentially equals to the Reidemeister torsion ${\rm Tor}(M,\rho)$. 
Lescop in \cite{Lescop} proved the equivalence for $b_1(M)=1$ and an abelian representation.
The author proved in \cite{Shimizu} that for $b_1(M)\ge 1$.
T.~Kitano and the author in \cite{KS} proved that both $d$ and ${\rm Tor}$ satisfy essentially same gluing formulas.
It is an evidence of essential equivalence between $d$ and ${\rm Tor}$.
We note that invariants, including the Chern-Simons perturbation theory and of course the defect $d$, whose origin are in the Chern-Simons quantum field theory should related to the Reidemeister torsion due to these physical backgrounds. See \cite{Witten}, \cite{Schwarz} for more details.

Both $d$ and ${\rm Tor}$ can be computed from combinatorial information of a Morse function and the representation. 
Take a Morse function on $M$.
We define a labeled graph called {\it Dehn graph} by extracting
information of the manifold by using the Morse function. 
Roughly speaking, a vertex of the Dehn graph corresponds to
a critical point of the Morse function, an edge is a trajectory connecting critical points.
We can compute both $d(M,\rho)$ and ${\rm Tor}(M,\rho)$ from the Dehn graph.

Let $K$ be a knot in $S^3$ and we denote by $E(K)$ the knot exterior of $K$.
In this article, we give a Dehn graph $\Gamma_K$ of $E(K)$ for a certain Morse function. 
We note that this Dehn graph is deeply related to the Dehn representation of the knot group. 
The Dehn graph $\Gamma_K$ is given by a combinatorial computation from a knot diagram of $K$. Hence finally we can compute both $d$ and ${\rm Tor}$ diagrammatically.


The organization of this paper is as follows.
In Section 2 we define a Dehn graph.
In Section 3 we describe the computation of several invariants including the defect $d$ by using the Dehn graph.
In Section 4 we give a diagrammatic construction of a Dehn graph for a certain Morse function on a knot exterior. 
In Section 5 we give an example of the diagrammatic method described in Sections 3 and 4. 
In Section 6 we give the Morse function used in Section 4. 
\subsection*{Acknowledgments}
The author expresses his
appreciation to Professor Seiichi Kamada for his valuable comments, in particular on the Dehn representation. 
The author would like to thank Professor Tetsuya Ito and Professor Ayumu Inoue for their helpful
comments, in particular on the diagrams used in Section 6.  
This work was partly supported by JSPS KAKENHI Grant Number JP18K13408.   

\section{Dehn graph}
Let $M$ be a compact oriented manifold with a base point $\infty\in M$.

Let $f:M\to \RR$ be a Morse function. If $M$ has a boundary $\partial M$,
we assume that $f^{-1}(0)=\partial M$.
Let $\grad f$ be a gradient like vector field satisfying
the Morse-Smale conditions.
Let $\Crit_i(f)$ denote the set of critical points of the Morse index $i$.  
Set $\Crit(f)=\sqcup_i\Crit_i(f)$.
For each $p\in \Crit(f)$, we denote its Morse index by $\ind(p)$.

For each critical point $p\in\Crit(f)$, we take a path $c_p$ from $p$ to $\infty$. 
\begin{defini}
A {\it bouquet} of $f$ is a collection $c_f=\{c_p\}_{p\in\Crit(f)}$ of such $c_p$.
\end{defini}
\includegraphics[width=10cm]{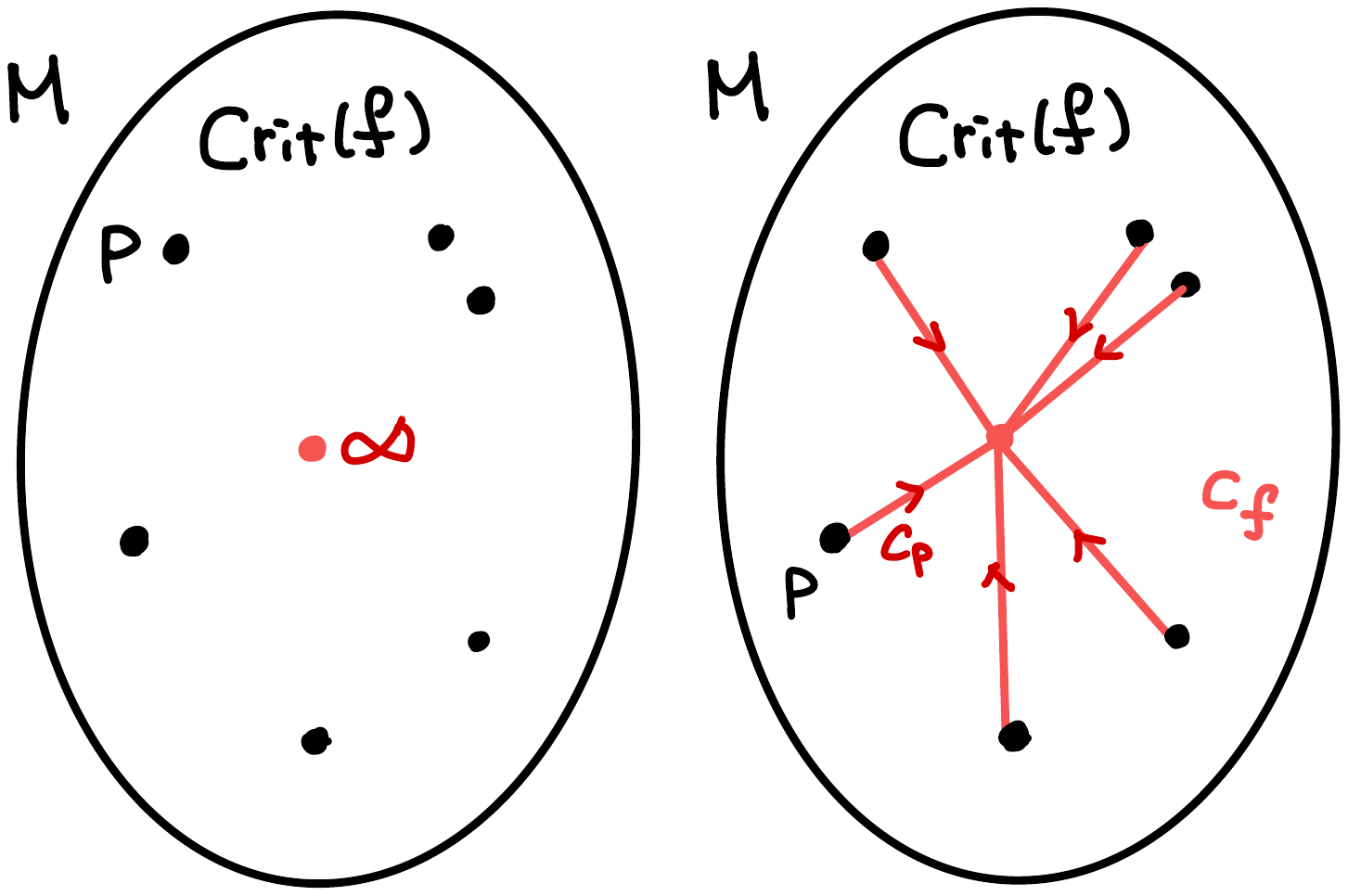}\label{pict1} 

A {\it trajectory} from $p\in\Crit(f)$ to $q\in\Crit(f)$ is an integral curve of $-\grad f$ from $p$ to $q$. We consider a trajectory as a path from $p$ to $q$, namely it has the downward direction with respect to $f$.
For $p,q\in \Crit(f)$ with $\ind(p)=\ind(q)+1$, let us denote by $\mathcal M(p,q)$ the set of trajectories from $p$ to $q$.
Let $\mathcal A_p$ the ascending disk of $p\in \Crit(f)$, namely the set of points on integral curves from $p$.
Let $\mathcal D_q$ the descending disk of $q\in \Crit(f)$ the set of points on integral curves converging to $q$.
We give orientation of $\mathcal A_p$ and $\mathcal D_p$ such that the condition $T_p\mathcal A_p\oplus T_p\mathcal D_p\cong T_pM$ holds for each $p\in {\rm Crit}(f)$. 
A trajectory $\gamma\in \mathcal M(p,q)$ is a part of $\mathcal A_p\cap \mathcal D_q$ as a manifold. Thus $\gamma$ has an orientation as a part of the intersection. 
Set $\ep(\gamma)=1$ if the orientation is from $p$ to $q$, otherwise
$\ep(\gamma)=-1$.

In this article, a {\it graph} is a 1-dimensional finite CW complex.
A $0$-simplex is called {\it vertex} and a $1$-simplex is called a {\it edge}.
For a graph $\Gamma$, we denote by $V(\Gamma)$ the set of vertexes. 
For vertexes $p,q\in V(\Gamma)$, we denote by $E(\Gamma)(p,q)$ the set of edges connecting $p$ and $q$. 
We denote by $E(\Gamma)=\sqcup_{p,q\in V(\Gamma)}E(p,q)$ the set of all edges.

We define a labeled graph $\Gamma_{f,c_f}$, called {\it Dehn graph}, as follows:
\begin{defini}[Dehn graph]
A {\it Dehn graph} $\Gamma=\Gamma_{f,c_f}$ of a Morse function $f$ and a bouquet $c_f$ is a labeled graph consisting of the following data:
\begin{itemize}
\item[(1)] $V(\Gamma)=\Crit(f)$.
\item[(2)] $E(\Gamma)(p,q)=\mathcal M(p,q)$ when $|\ind(p)-\ind(q)|=1$.
Otherwise, $E(\Gamma)(p,q)=\emptyset$.
\item[(3)] Each $p\in V(\Gamma)$ has a label $l(p)=\ind(p)$.
\item[(4)] Each $\gamma\in \mathcal M(p,q)=E(\Gamma)(p,q)$ has a label 
$l(\gamma)=\ep(\gamma)[c_q\circ \gamma\circ c_p^{-1}]\in \ZZ[\pi_1(M,\infty)]$. 
Here $[c_q\circ \gamma\circ c_p^{-1}]\in \pi_1(M,\infty)$ is a class represented by a loop $c_q\circ \gamma\circ c_p^{-1}$.
\end{itemize}  
\end{defini}


\section{Invariants computed from a Dehn graph} 
Let $M$ be a compact oriented manifold.
Let $\rho:\pi_1(M)\to GL_F(V)$ be a finite dimensional linear representation, where $V$ is a vector space over a field $F$.
Let $V_{\rho}$ denotes the corresponding local system.
Take a Morse function $f:(M,\partial M)\to (\RR_{\le 0},0)$ and a bouquet $c_f$.
Let $\Gamma$ be a Dehn graph of $f$ and $c_f$.
In this section, we compute some invariants (the Morse-Smale complex, the Reidemeister torsion and the defect) of a pair $(M,\rho)$ from the Dehn graph $\Gamma$.

\subsection{Morse-Smale complex $C_*^{f,\rho}$ and the homology groups $H_*(M;V_{\rho})$}
One can compute the Morse-Smale complex 
$$C^{f,\rho}_*=\cdots \to C_{i+1}\stackrel{\partial_{i+1}}{\to} C_i\cdots$$
of $(M,\rho)$ from $\Gamma$ as follows:
We assign $V_p$, which is a copy of $V$, to each vertex $p\in V(\Gamma)$.
The complex $C_i$ is the direct sum of $V_p$ of the vertex $p$ with $l(p)=\ind(p)=k$: 
$$C_k=\oplus_{p\in V(\Gamma)_k}V_p, ~~ V_{p}=V.$$
Here $V(\Gamma)_k=\{p\in V(\Gamma)\mid l(p)=k\}$.
Set
$$\partial_{p,q}=\sum_{\gamma \in E(\Gamma)(p,q)}
\rho(l(\gamma)):V_{p}\to V_{q}.$$
Here $\rho(-a)$ for $a\in \pi_1(M,\infty)$ is defined as $-\rho(a)$. 
Then the boundary operator $\partial_k:C_k\to C_{k-1}$ is written as   
$$\partial_k=\sum_{p\in V(\Gamma)_k,q\in V(\Gamma)_{k-1}}\partial_{p,q}.$$
Of course, by the Morse theory,
$$H_*(M;V_{\rho})=H_*(C^{f,\rho}_*).$$
\vskip4mm
\subsection{Reidemeister torsion ${\rm Tor}(M,\rho)\in F^{\times}/
\pm\det\rho(\pi_1(M))$}
In this section, we review the computation of the Reidemeister torsion 
${\rm Tor}(M,\rho)$ from the Morse-Smale complex $C^{f,\rho}_*$. 
In this and the next section, we assume that $H_*(M;V_{\rho})=0$, namely $C^{f,\rho}_*=\cdots\to C_{i+1}\to C_{i}\to \cdots$ is an exact sequence.

Take a family of homomorphisms $G=\{G_i:C_{i-1}\to C_i\}_{i}$ satisfying
$$\partial_{i+1}\circ G_{i+1}+G_i\circ \partial_i=1_{C_i}$$
for any $i$.
$G$ is called a {\it combinatorial propagator}. 

Let $C_{\rm even}=\oplus_kC_{2k}, C_{\rm odd}=\oplus_kC_{2k+1}$.
Then $\partial+G$ is an isomorphism:
$$\partial+G:C_{\rm even}\to C_{\rm odd}.$$
We fix a basis of $V$ over $F$. It gives basis of $C_{\rm odd}$ and $C_{\rm even}$.
Then the determinant $\det(\partial+G)$ of $\partial+G$ takes a value in $F^{\times}=F\setminus\{0\}$.
There are ambiguities of the choice of $G$ and the choice of a basis of $V$. However, it is easily checked that
$\det(\partial+G)$ is independent from these choices.
Furthermore, it is known that $\det(\partial+G)$ modulo $\pm\det\rho(\pi_1(M))$ is independent of the choice of the Morse function $f$ and the bouquet $c_f$.  
\begin{defini}[Reidemeister torsion]
${\ Tor}(M,\rho)=\det(\partial +G)\in F^{\times}/\pm\det\rho(\pi_1(M)).$
\end{defini}

\subsection{Defect $d(M,\rho)\in H_1(M;V_{\rho^*}\otimes V_{\rho})/H_1(M;\ZZ)$}
Both a Dehn graph $\Gamma$ and a bouquet $c_f$ can be considered as a singular 1-chain. In fact, $\Gamma$ combined with $c_f$ with appropriate local coefficients determined from the labels of edges gives a 1-cycle.   
Then the invariant $d(M,\rho)$, which is called {\it defect}, is given as a homology class represented by the cycle.
$$d(M,\rho)\in H_1(M;V_{\rho^*}\otimes V_{\rho})/H_1(M;\ZZ).$$
Here $\rho^*=\null^t\rho^{-1}$ is the dual representation of $\rho$.
See \cite[Section 4]{Shimizu} for more details of $d$.

For a singular 1-simplex $c:[0,1]\to M$ and 
$\alpha\in V_{\rho^*}\otimes V_{\rho}$, we have a 1-chain
$$c\alpha\in C_1(M;V_{\rho^*}\otimes V_{\rho}).$$
In particular, for $[c]\in\pi_1(M;\infty)$
and $\alpha \in V_{\rho^*}\otimes V_{\rho}$, a 1-chain 
$[c]\alpha$ makes sense.   

Take a combinatorial propagator $G=\{G_i:C_{i-1}\to C_i\}_i$.
For $p,q\in V(\Gamma)$ with $l(p)={\rm ind}(p)=i$ and $l(q)=i-1$, let $G_{p,q}=\pi_{p}\circ G_i|_{V_q}:V_q\to V_p$ be a part of $G$, where 
$\pi_p:C_{i}\to V_p$ is the projection. 
Then for each edge $\gamma\in E(\Gamma)(p,q)$, 
the composition $\rho(l(\gamma))\circ G_{p,q}$ is an element of 
$V_{\rho^*}\otimes V_{\rho}$.
Let $\rho(l(\gamma))\circ G$ denote $\rho(l(\gamma))\circ G_{p,q}$ shortly.
Then $d(M,\rho)$ is described as follows (\cite[Proposition 6.8]{Shimizu}):
$$d(M,\rho)=\left[\sum_{\gamma\in E(\Gamma)}
l(\gamma)(\rho(l(\gamma))\circ G)\right]\in H_1(M;V_{\rho^*}\otimes V_{\rho})/H_1(M;\ZZ).$$
For any $\gamma$ with $l(\gamma)=\pm1$, 
$l(\gamma)\rho((l(\gamma))\circ G)=0$. Then we have
$$d(M,\rho)=\left[\sum_{\gamma\in E(\Gamma), l(\gamma)\not=\pm 1}
l(\gamma)(\rho(l(\gamma))\circ G)\right].$$

\section{A diagrammatic method for a Dehn graph of a knot exterior}\label{Sect:combi}
Let $K\subset \RR^3\subset S^3=\RR^3\cup\{\infty\}$ be a knot.
Let $E(K)$ denote the knot exterior, namely $E(K)=S^3\setminus N(K)$, where $N(K)$ is an open tubular neighborhood of $K$.
In this section, we give a combinatorial method to construct a Dehn graph $\Gamma_K$ of a certain Morse function on $E(K)$ and a certain bouquet of the Morse function.
We give the Morse function and the bouquet in Section~\ref{Sect:Morse}.

Let $\pi:\RR^3\to \RR^2, (x,y,z)\mapsto (x,y)$ be the projection.
Let $D\subset \RR^2$ be a regular knot diagram of $K$ with the crossings $P_1,\ldots,P_k$.
$\RR^2\setminus \pi(K)$ consists of one unbounded region and $(k+1)$-bounded regions.
Let $Q_{\infty}$ be the unbounded region and $Q_1,\ldots,Q_{k+1}$ be the bounded regions.
The knot diagram $D$ has $k$ arcs.
Let $a,b,c,\ldots$ be these arcs.
Then a Wirtinger representation of $\pi_1(E(K))$ is given as
$$\pi_1(E(K))=\langle a,b,c,\ldots\mid {\rm Rel}\rangle,$$
where ${\rm Rel}$ is a set of relations corresponding to the crossings.

\includegraphics[width=8cm]{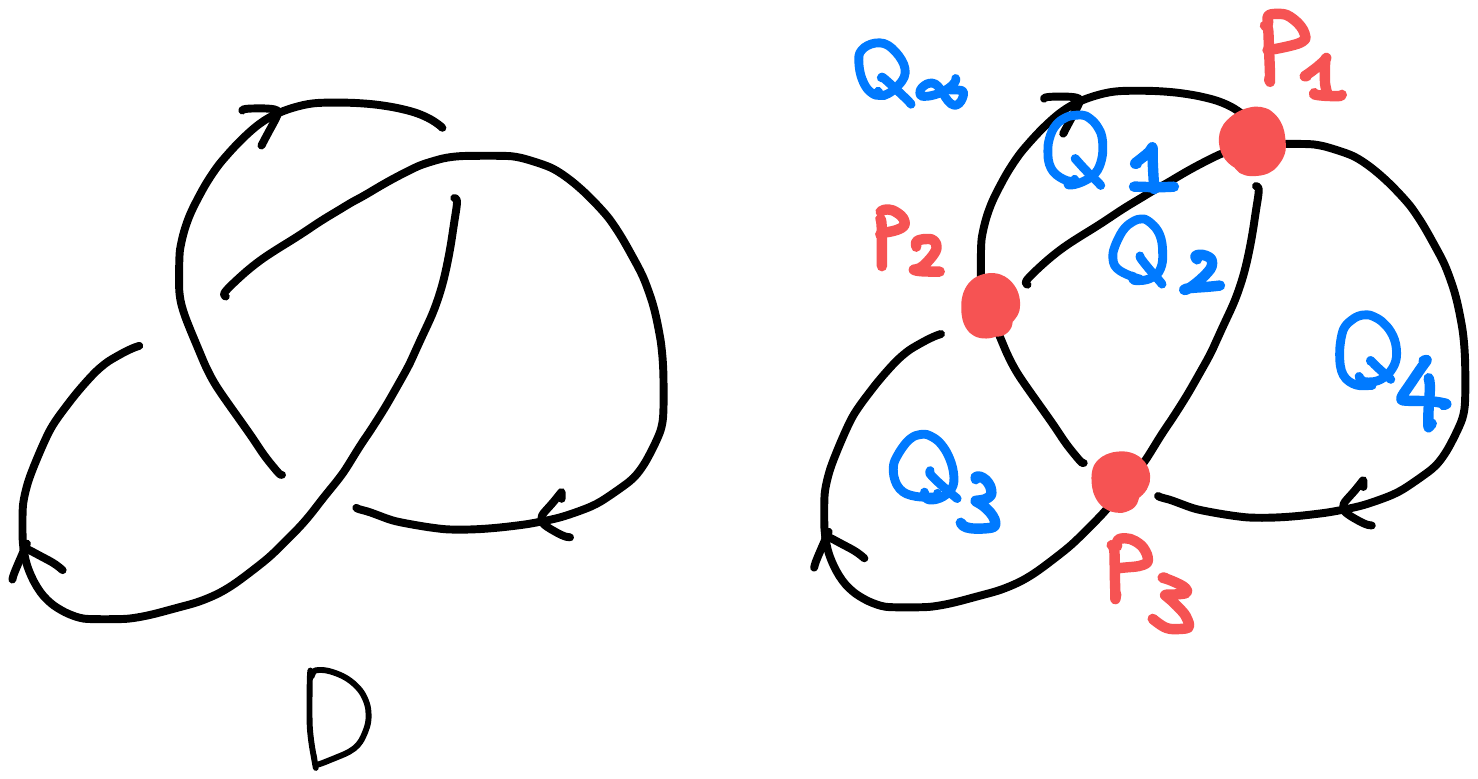}
\includegraphics[width=4cm]{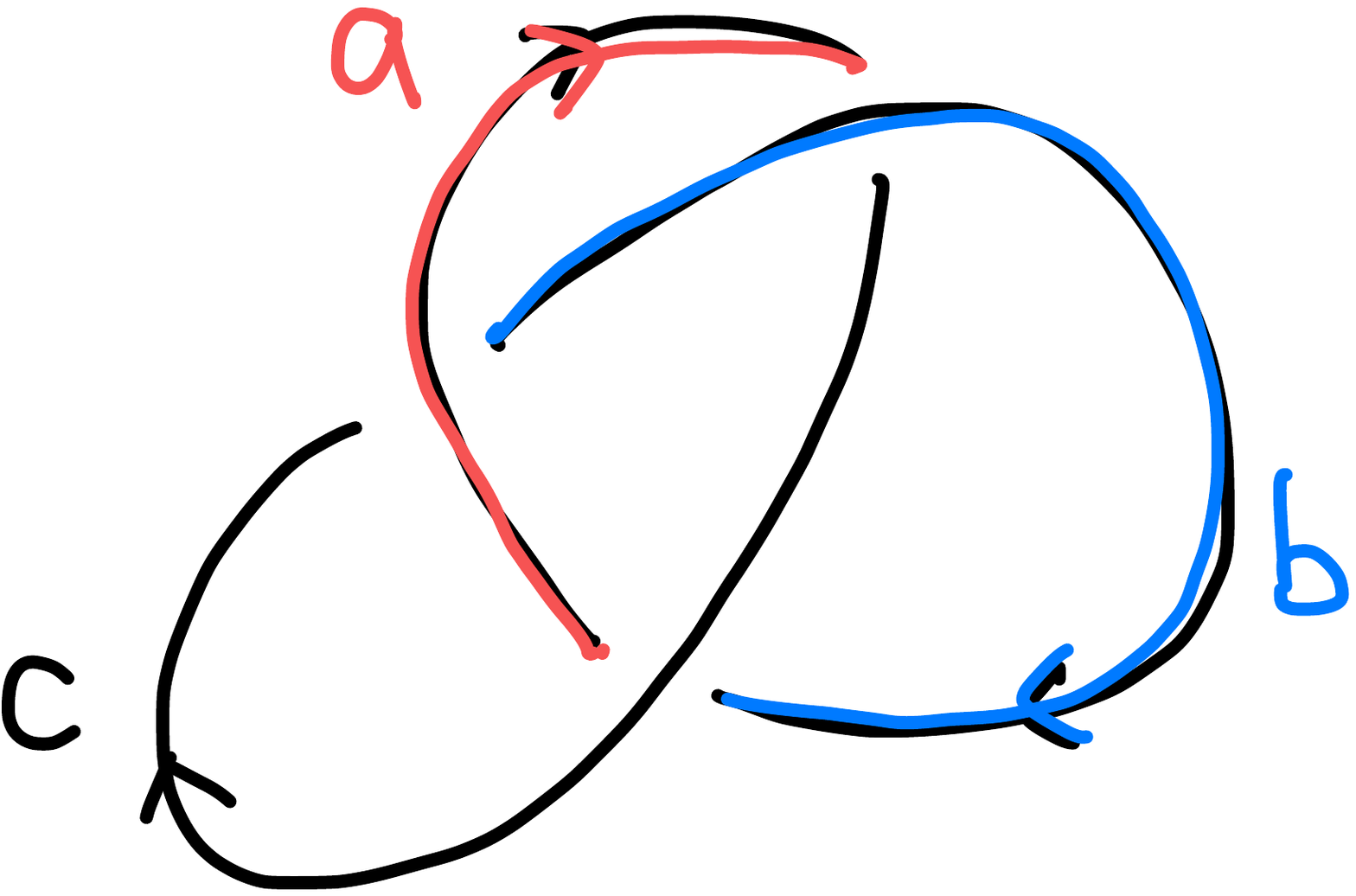}

To compute a Dehn graph, we introduce two diagrams $D_1,D_2$ given from $D$.
\subsection{The diagram $D_1$}\label{D_1}
For each crossing $P_i$, there are four corners $R_i^1,\ldots, R_i^4$ of regions near $P_i$.  
We assign $n(R_i^1),\ldots, n(R_i^4)\in \ZZ[\pi_1(E(K))] $ for each corner by the following rule (R):
\begin{itemize}
\item[(R)] Let $x$ be the over arc at $P_i$.
Before the crossing $P_i$, we assign $-x$ and $1$ to
the corner on the left side and the right side of the over arc $x$ respectively.
After the crossing $P_i$, we assign $x$ and $-1$ to the corner on the left side and the right side of $x$ respectively.
\end{itemize}

\includegraphics[width=8cm]{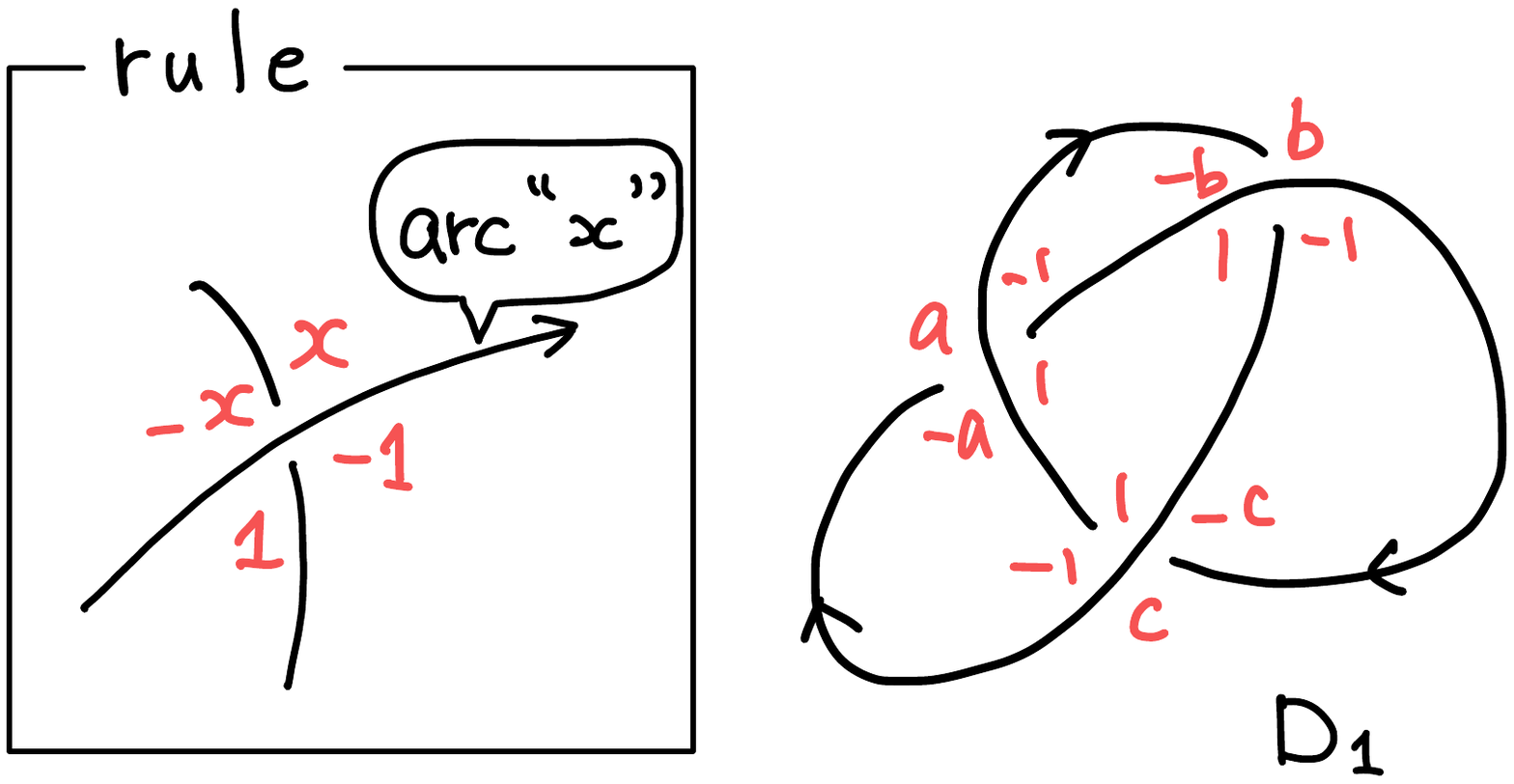}

\noindent
{\bf Remark.}~
The diagram $D_1$ is essentially given by M. Cohen, O. T. Dasbach and H. M. Russell in \cite{Cohen} as a twisted Alexander graph.
The diagram $D_1$ can be considered as a generalization of
the Alexander diagram for a diagrammatic computation of the Alexander polynomial from the Dehn representation of the knot group.

\vskip3mm
\subsection{The diagram $D_2$}\label{D_2}
The diagram 2 is given as follows.
We assign an element $l(Q_i)$ of $\pi_1(E(K))$ for each region $Q_i$ for $i=1,\ldots,k+1,\infty$ by the following rules (R1) and (R2):
\begin{itemize}
\item[(R1)] $l(Q_{\infty})=1$.
\item[(R2)] If an arc $x$ divides $Q_{i}$ and $Q_j$ so that $Q_i$ is on the left side of the arc, then $l(Q_j)=xl(Q_i)$.
\end{itemize}
\includegraphics[width=8cm]{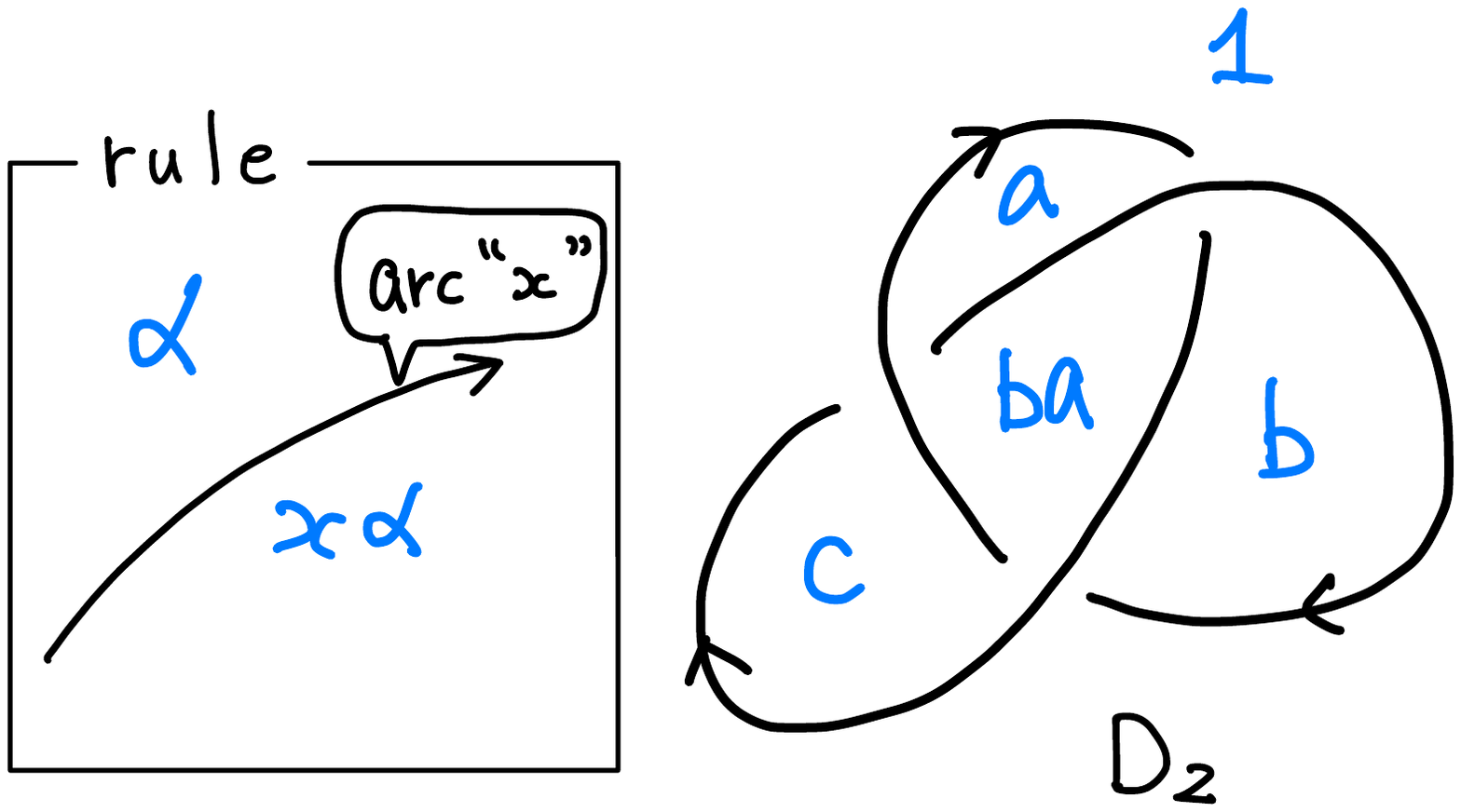}
\subsection{A Dehn graph $\Gamma_K$ computed from $D_1$ and $D_2$}
We give a Dehn graph ${\Gamma}_K$ of a certain Morse function and a certain bouquet by using the diagrams $D_1$ and $D_2$.
\begin{itemize}
\item $\Gamma_K$ has $2k+2$-vertexes $p_1,\ldots,p_{k}, q_1,\ldots,q_{k+1}$ and $\infty$.  
$$V(\Gamma)=\{p_1,\ldots, p_k,q_1,\ldots,q_{k+1},\infty\}.$$ 
The label of each $p_1,\ldots, p_k$ is $2$, that of $q_1,\ldots, q_{k+1}$ are $1$ and that of $\infty$ is $0$.
The vertex $p_i$ is corresponding to the crossing $P_i$.
The vertex $q_j$ is corresponding to the region $Q_j$. 
\item If the crossing $P_i$ is on the boundary of the region $Q_j$, we put an edge $\gamma_{p_i,q_j}$ connecting the vertex $p_i$ to the vertex $q_j$:
$$E(\Gamma_K)(p_i,q_j)=\{\gamma_{p_i,q_j}\}.$$ 
The label $l(\gamma_{p_i,q_j})$ is the label assigned to the corner between $P_i$ and $Q_j$ in the diagram $D_1$.
\item We put two edges connecting each $q_j$ and $\infty$:
$$E(\Gamma_K)(q_j,\infty)=\{\gamma_{q_j,\infty}^+,\gamma_{q_j,\infty}^-\}.$$
The label of $\gamma_{q_j,\infty}^+$ is $1$.
The label of $\gamma_{q_j,\infty}^-$ is $(-1)$ times
the label of the region $Q_j$ given in the diagram $D_2$.  
\item There are no edges other than the above.
\end{itemize}
\includegraphics[width=8cm]{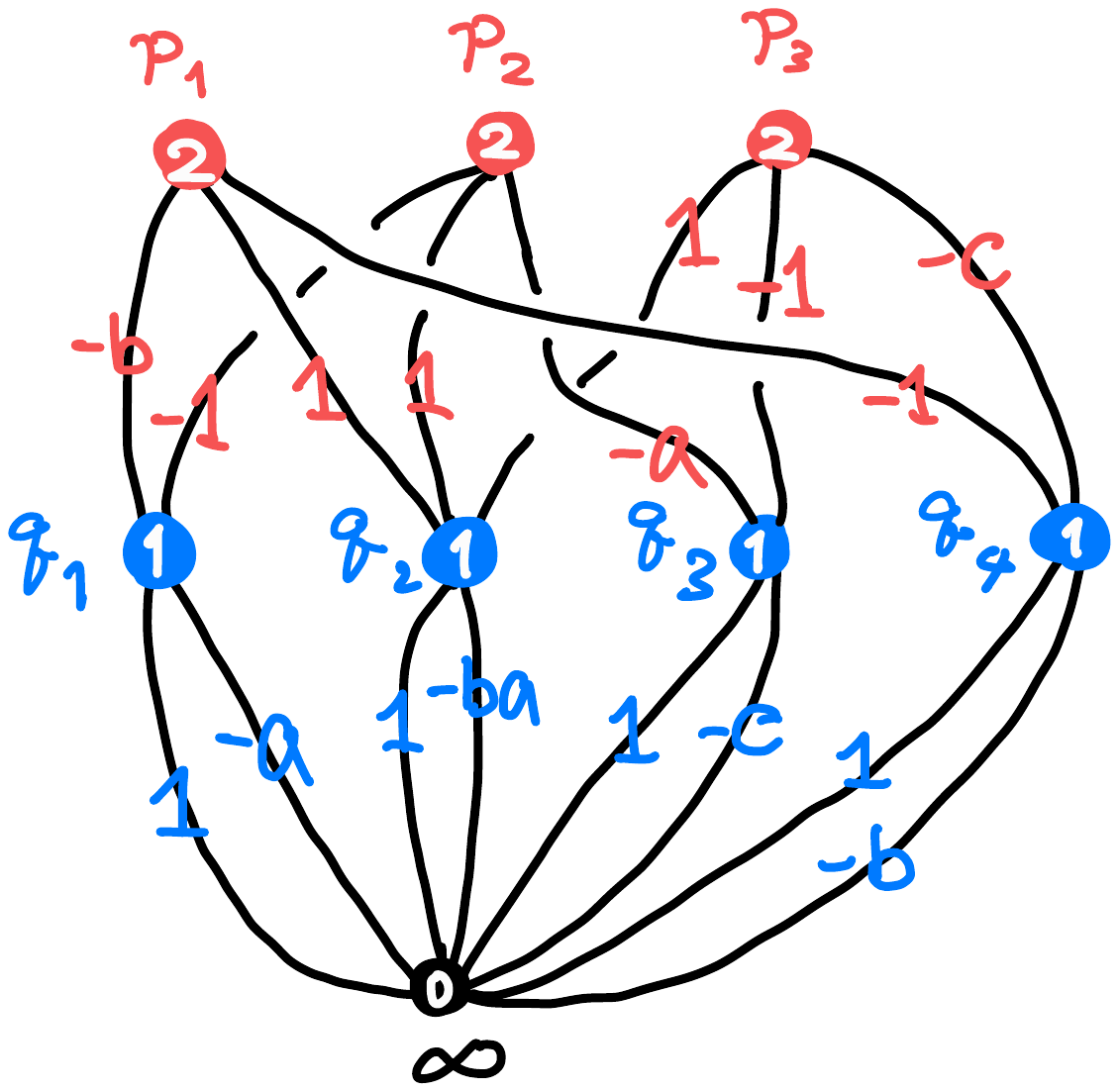}\label{pict6}


\section{Example} 
As an example, we compute ${\rm Tor}(E(K),\rho_{0})$ and $d(E(K),\rho_{0})$
for the trefoil knot $K$ with the maximal abelian representation 
$$\rho_{0}:\pi_1(E(K))\to H_1(E(K);\ZZ)=\langle t\rangle \ni t\mapsto t\in Q(t)^{\times}=GL(Q(H)).$$ 
Here $Q(t)$ is the quotient field of the polynomial ring $\QQ[t]$.
We consider that $Q(H)$ is a 1-dimensional vector space over a field $Q(H)$.
We note that both ${\rm Tor}(E(K),\rho_{0})$ and $d(E(K),\rho_0)$ can be computed from Alexander polynomial.
This example is just an example for our combinatorial method.

We use diagrams $D$, $D_1$ and $D_2$ described in the previous section.

Then we have the Dehn graph $\Gamma_K$ described in the previous section.
\subsection{the Morse-Smale complex $C_*^{f,\rho_0}$}\label{51}
We first compute the Morse-Smale complex $C_*^{f,\rho_0}
=0\to C_2\stackrel{\partial_2}{\to} C_1\stackrel{\partial_1}{\to}
C_0\to 0$.
\begin{eqnarray*}
C_2&=&Q(t)_{p_1}\oplus Q(t)_{p_2}\oplus Q(t)_{p_3}\cong Q(t)^{\oplus 3},\\
C_1&=&Q(t)_{q_1}\oplus Q(t)_{q_2}\oplus Q(t)_{q_3}\oplus Q(t)_{q_4}\cong Q(t)^{\oplus 4},\\
C_0&=&Q(t)_{\infty}\cong Q(t).
\end{eqnarray*}
We denote by $p_1=(1,0,0)\in Q(t)^{\oplus 3}=C_2$. 
We use notations $p_2,p_3,q_1,\ldots,q_4$ and $\infty$ in the same manner. 
Then $\{p_1,p_2,p_3\}$, $\{q_1,\ldots, q_4\}$ and $\{\infty\}$ are basis of $C_2$, $C_1$ and $C_0$ respectively. 
The representation $\rho_0$ satisfies
$\rho_0(a)=\rho_0(b)=\rho_0(c)=t$.
Therefore, by using these basis, the boudary homomorphisms can be written as follows:
\[ \partial_2(p_1,p_2,p_3)=(q_1,q_2,q_3,q_4)
\left(\begin{array}{ccc}
-t&-1&0\\
1&1&1\\
0&-t&-1\\
-1&0&-t
\end{array}
\right),\]
\[ \partial_1(q_1,q_2,q_3,q_4)=(\infty)
\left(\begin{array}{cccc}
1-t&1-t^2&1-t&1-t
\end{array}
\right).\]
\subsection{A propagator $G$}
It is easily checked that $\{q_1,\partial_2(p_1),\partial_2(p_2),\partial_1(p_3)\}$ is a basis of $C_2$.
Similarly, $\{\partial_1(q_1)\}$ is a basis of $C_0$.
Let $G_2:C_1\to C_2$ and $G_1:C_0\to C_1$ be
\[G_2(q_1,\partial_2(p_1),\partial_2(p_2),\partial_2(p_3))
=(p_1,p_2,p_3)\left(
\begin{array}{cccc}
0&1&0&0\\
0&0&1&0\\
0&0&0&1
\end{array}
\right) \text{and}\]
\[G_1(\partial_1(q_1))=(q_1,\partial_2(p_1),\partial_2(p_2),\partial_2(p_3))
\left(
\begin{array}{c}
1\\
0\\
0\\
0
\end{array}
\right).\]
These $G_2,G_1$ satisfy
$G_2\circ\partial_2={\rm id}_{C_2},
\partial_2\circ G_2+G_1\circ\partial_1={\rm id}_{C_1}$
and $\partial_1\circ G_1={\rm id}_{C_0}$.
Namely $G=\{G_2,G_1\}$ is a propagator.
For the convenience of the computations, we give a description of $G$ under the basis given in Section~\ref{51}.
Since $$(q_1,\partial_2(p_1),\partial_2(p_2),\partial_2(p_3))
=(q_1,q_2,q_3,q_4)\left(
\begin{array}{cccc}
1&-t&-1&0\\
0&1&1&1\\
0&0&-t&-1\\
0&-1&0&-t
\end{array}
\right),$$
we have
\begin{eqnarray*}
G_2(q_1,q_2,q_3,q_4)&=&(p_1,p_2,p_3)
\left(
\begin{array}{cccc}
0&1&0&0\\
0&0&1&0\\
0&0&0&1
\end{array}
\right)
\left(\begin{array}{cccc}
1&-t&-1&0\\
0&1&1&1\\
0&0&-t&-1\\
0&-1&0&-t
\end{array}\right)^{-1}\\&=&(p_1,p_2,p_3)\frac{1}{t^2-t+1}
\left(
\begin{array}{cccc}
0&t^2&t&t-1\\
0&1&1-t&1\\
0&-t&-1&-t
\end{array}
\right) \text{and}
\end{eqnarray*}
\begin{eqnarray*}
G_1(\infty)&=&(q_1,q_2,q_3,q_4)
\left(\begin{array}{cccc}
1&-t&-1&0\\
0&1&1&1\\
0&0&-t&-1\\
0&-1&0&-t
\end{array}\right)
\left(
\begin{array}{c}
\frac{1}{1-t}\\
0\\
0\\
0
\end{array}\right)\\
&=&(q_1,q_2,q_3,q_4)
\left(
\begin{array}{c}
\frac{1}{1-t}\\
0\\
0\\
0
\end{array}\right).
\end{eqnarray*}

\subsection{Reidemeister torsion ${\rm Tor}(E(K),\rho_0)$}
An isomorphism $\partial +G=\partial_2+G_1:C_{\rm even}\to C_{\rm odd}$ is described as follows:
$$(\partial_2+G_1)(p_1,p_2,p_3,\infty)
=(q_1,q_2,q_3,q_4)\left(
\begin{array}{cccc}
-t&-1&0&\frac{1}{t-1}\\
1&1&1&0\\
0&-t&-1&0\\
-1&0&-t&0
\end{array}
\right).
$$
Therefore the Reidemeister torsion ${\rm Tor}(E(K),\rho_0)$ is computed as
$${\rm Tor}(E(K),\rho_0)=\det
\left(
\begin{array}{cccc}
-t&-1&0&\frac{1}{t-1}\\
1&1&1&0\\
0&-t&-1&0\\
-1&0&-t&0
\end{array}
\right)=\frac{1}{1-t}(t^2-t+1).$$
\subsection{Defect $d(E(K),\rho_0)$}
Note that $H_1(E(K);Q(t)_{\rho_0^*}\otimes Q(t)_{\rho_0})\cong
H_1(E(K);\ZZ)\otimes Q(t)\cong Q(t)$. Under this isomorphism,
$H_1(E(K);\ZZ)\subset H_1(E(K);Q(t)_{\rho_0^*}\otimes Q(t)_{\rho_0})$
corresponds to $\ZZ \subset Q(t)$. Then
$$d(E(K),\rho_0)\in H_1(E(K);Q(t)_{\rho_0}^*\otimes Q(t)_{\rho_0})/
H_1(E(K);\ZZ)=Q(t)/\ZZ.$$
For the convenience for the computation, we give 
the following figure which shows a propagator corresponding to each edge.

\includegraphics[width=8cm]{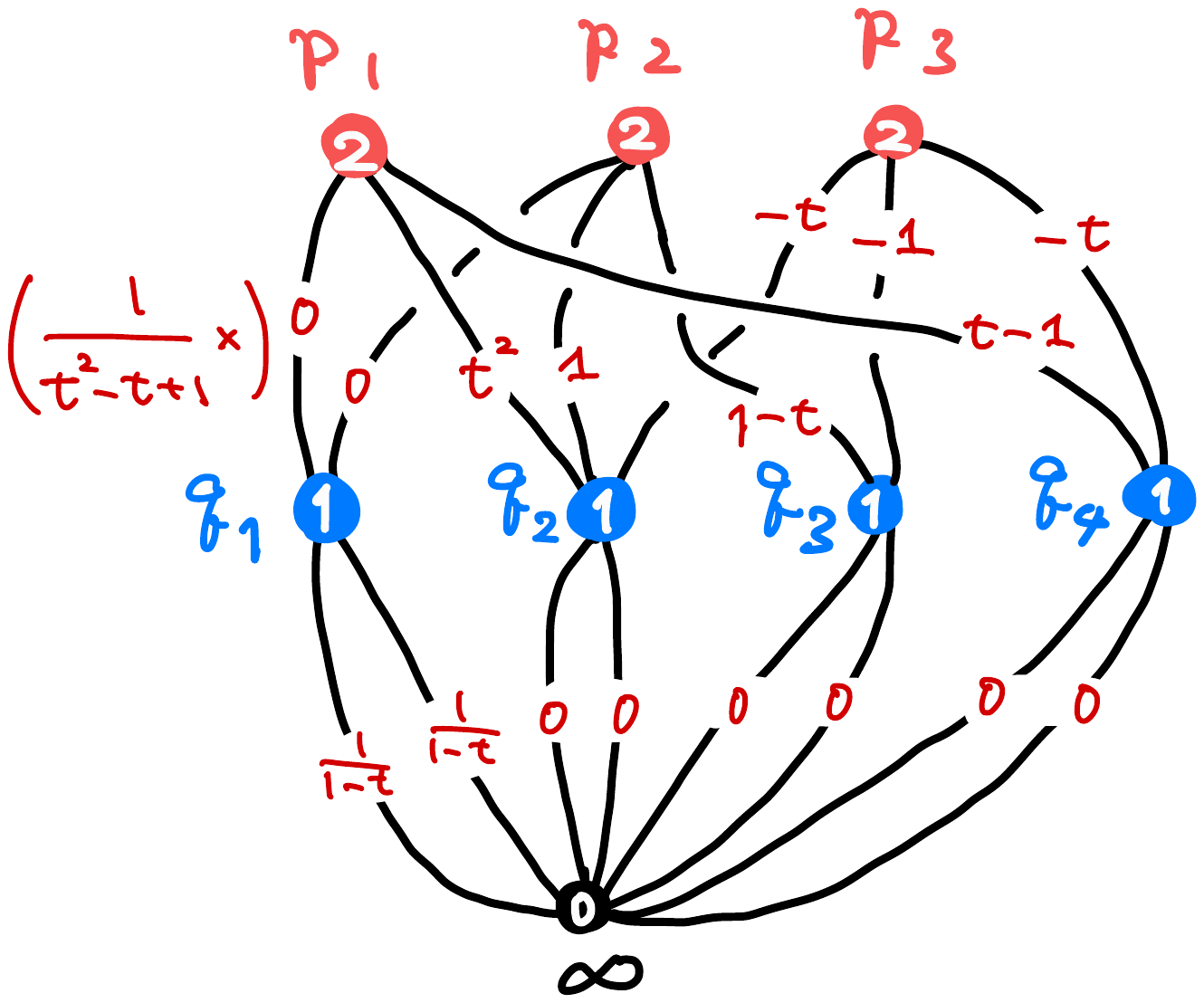}

\begin{eqnarray*}
d(E(K),\rho_0)&=&\left[\sum_{\gamma\in E(\Gamma), l(\gamma)\not=\pm1}
[l(\gamma)]\cdot \rho_0(l(\gamma))\circ G)\right]\\
&=&\left[ \frac{-t(1-t)}{t^2-t+1}+\frac{-t(-t)}{t^2-t+1}+\frac{-t}{1-t}\right]\\
&=& \left[\frac{2t^2-t}{t^2-t+1}-\frac{t}{t-1}\right]\\
\Big(&=&t\frac{d}{dt}\log\left\{ \frac{1}{1-t}(t^2-t+1)\right\}=t\frac{d}{dt}\log {\rm Tor}(E(K),\rho_0)\Big). 
\end{eqnarray*}

\section{A Morse function on $E(K)$}\label{Sect:Morse}
The labeled graph $\Gamma_K$ described in Section~\ref{Sect:combi} is a Dehn graph of a certain Morse function $f:(E(K),\partial E(K))\to (\RR_{\le 0},0)$ on $E(K)$ and a bouquet $c_f$ for $f$. 
Namely $\Gamma_K=\Gamma_{f,c_f}$.
In this section, we give $f$ and $c_f$.
Then we give a proof of $\Gamma_K=\Gamma_{f,c_f}$.
\vskip3mm
\noindent
\subsection{A descriptions of $f$ and $c_f$}
\noindent
\fbox{Morse function}
\vskip2mm
Recall that our knot $K$ is in $\RR^3=S^3\setminus \{\infty\}$ and $\pi:\RR^3\to \RR^2=\RR^2\times \{0\}\subset \RR^3$ is the projection. 
$P_1,\ldots,P_k$ are the crossings of $\pi(K)\subset \RR^2$.
We assumed that $K$ is on $\RR^2$ except around the crossings and the preimage of each crossing is on both side of $\RR^2$. 

We take a Morse function $f:(E(K),\partial E(K))\to (\RR_{\le 0},0)$
 and a gradient like vector field $\grad f$ satisfying the following properties:
\begin{itemize}
\item $\Crit(f)=\Crit_2(f)\sqcup \Crit_1(f)\sqcup \Crit_0(f)$,
$\Crit_2(f)=\{P_1,\ldots, P_k\}\subset\RR^2$, 
$\Crit_1(f)=\{\widetilde{Q}_1,\ldots,\widetilde{Q}_{k+1}\}\subset\RR^2$ and
$\Crit_0(f)=\{\infty\}$. Here each critical point $\widetilde{Q}_j$ is around the center of the region $Q_j$. 
\item $f^{-1}(0)=\partial E(K)$.
\item $f(P_1)=\cdots=f(P_k)=-1$,
\item $f(\widetilde{Q}_1)=\cdots=f(\widetilde{Q}_{k+1})=-2$,
\item $f(\infty)=-3$.
\item For each crossing $P_i$, let $L_i$ be a line segment in $\RR^3$
connecting two preimages of $P_i$ under the projection form $K$ to $\pi(K)$. 
For any $0\ge t> -1$, $f^{-1}(t)\cap L_i$ consists of just two points.
As $t$ approaches $-1$, the two points $f^{-1}(t)\cap L$ get closer.
When $t=-1$, two points meet and then vanishing, namely $f^{-1}(t)\cap L=\emptyset$ for $-1>t$.
\item For $-1>t>-2$, $f^{-1}(t)$ is a boundary of a retract neighborhood of $K\cup (\sqcup_{i=1}^{k} L_i)$. 
In other word, for $-1>t>2$, $f^{-1}(t)$ is a boundary of a retract neighborhood of $\pi(K)$.

\includegraphics[width=9cm]{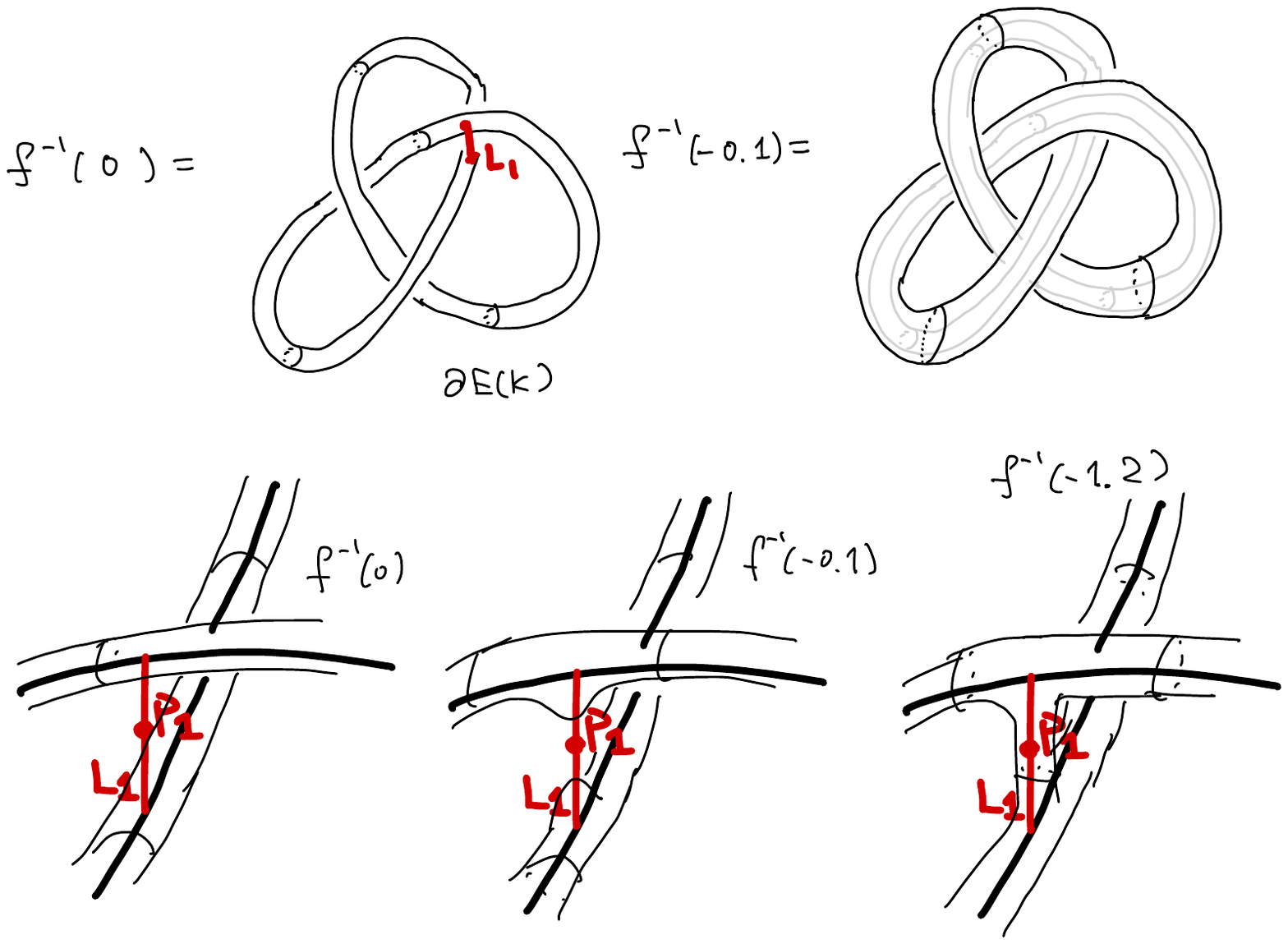}\label{pict7}
\includegraphics[width=7cm]{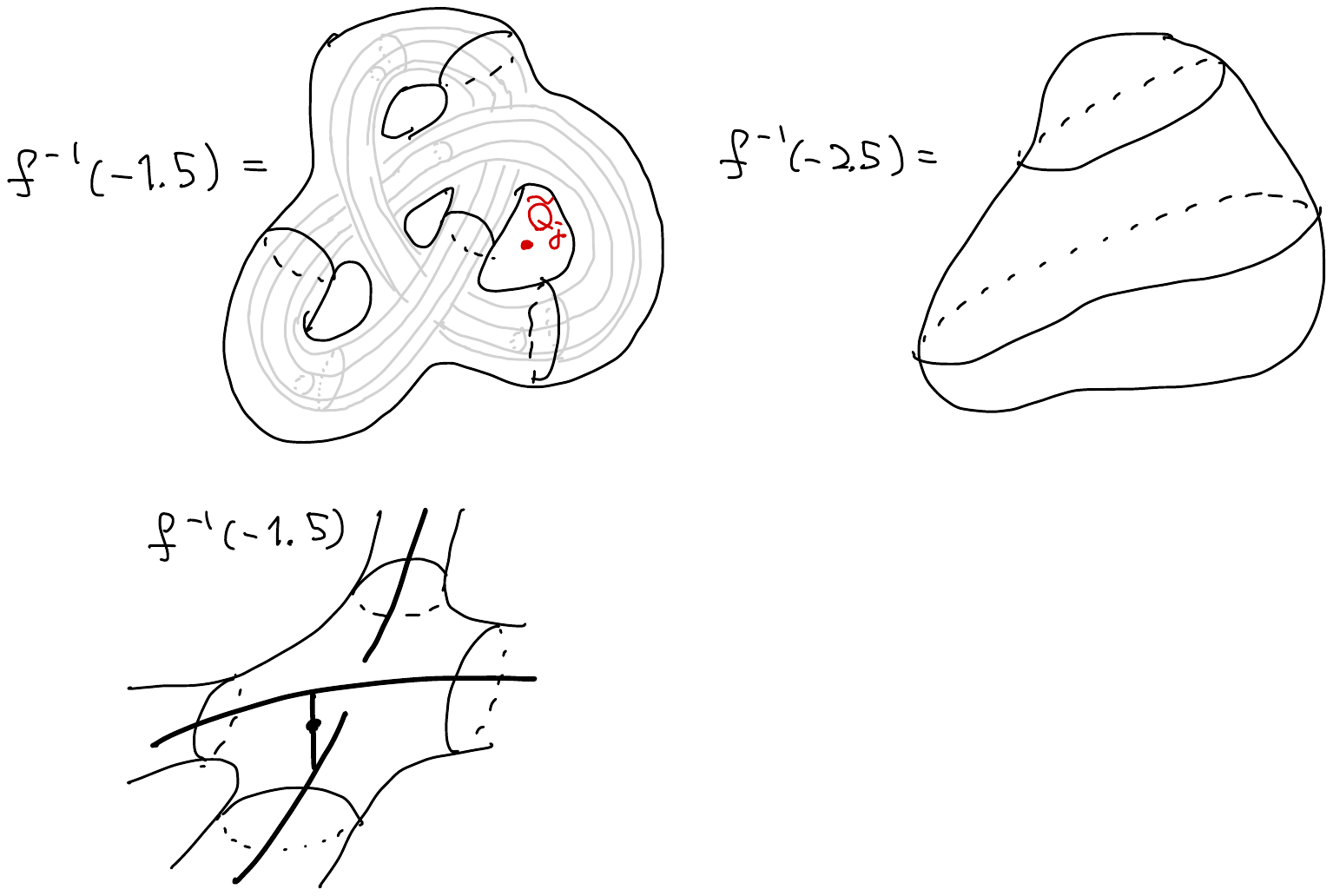}\label{pict72}

\item For $-1>t>-2$, $f^{-1}(t)$ is a surface of genus $k+1$.
The genus correspond to the critical points $\widetilde{Q}_1,\ldots, \widetilde{Q}_{k+1}$.
When $t=-2$, at the critical points $\widetilde{Q}_1,\ldots, \widetilde{Q}_{k+1}$, all genus are vanishing.
Then for $-2>t>-3$, $f^{-1}(t)$ is homeomorphic to the sphere $S^2$.

\item At $t=-3$, $f^{-1}(-3)=\{\infty\}$. 
\item $f^{-1}(t)=\emptyset$ for any $t<-3$. 
\end{itemize}
\vskip3mm
\noindent
\fbox{Bouquet $c_f$ for $f$}
\vskip3mm
Recall that a bouquet $c_f$ is a collection of paths $c_p$ from each critical point $p\in \Crit(f)$ to the base point $\infty$ in $E(K)$:
$c_f=\cup_{p\in\Crit(f)}c_p.$

To give a path $c_{P_i}$ for $i=1,\ldots,k$, we introduce a moving point $\alpha=\alpha(s), s\in [0,\infty]$, $\alpha(0)=P_i, \alpha(\infty)=\infty$.
The path $c_{P_i}$ is given as an orbit of $\alpha(s)$. 
After departing $p_i$, during $0\le s<1$, $\alpha$ moves slightly, in $\RR^2$, in the direction to the region up ahead on the right of the over arc as in the following picture. For $1\le s$, $\alpha(s)=\alpha(1)+s(0,0,1)$. Then $\alpha$ finally arrives at $\infty$.

\includegraphics[width=5cm]{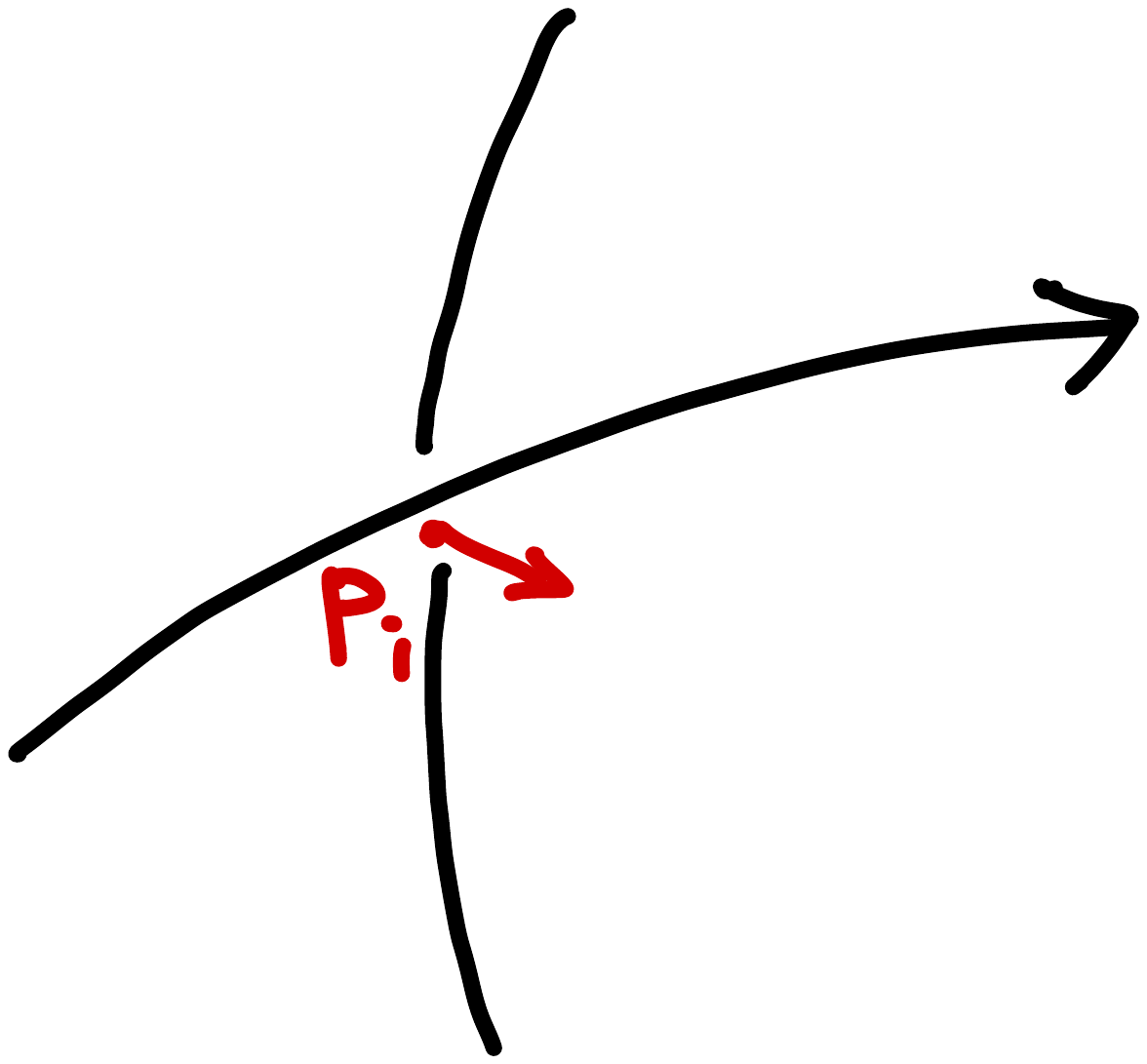}\label{pict91}
\includegraphics[width=8cm]{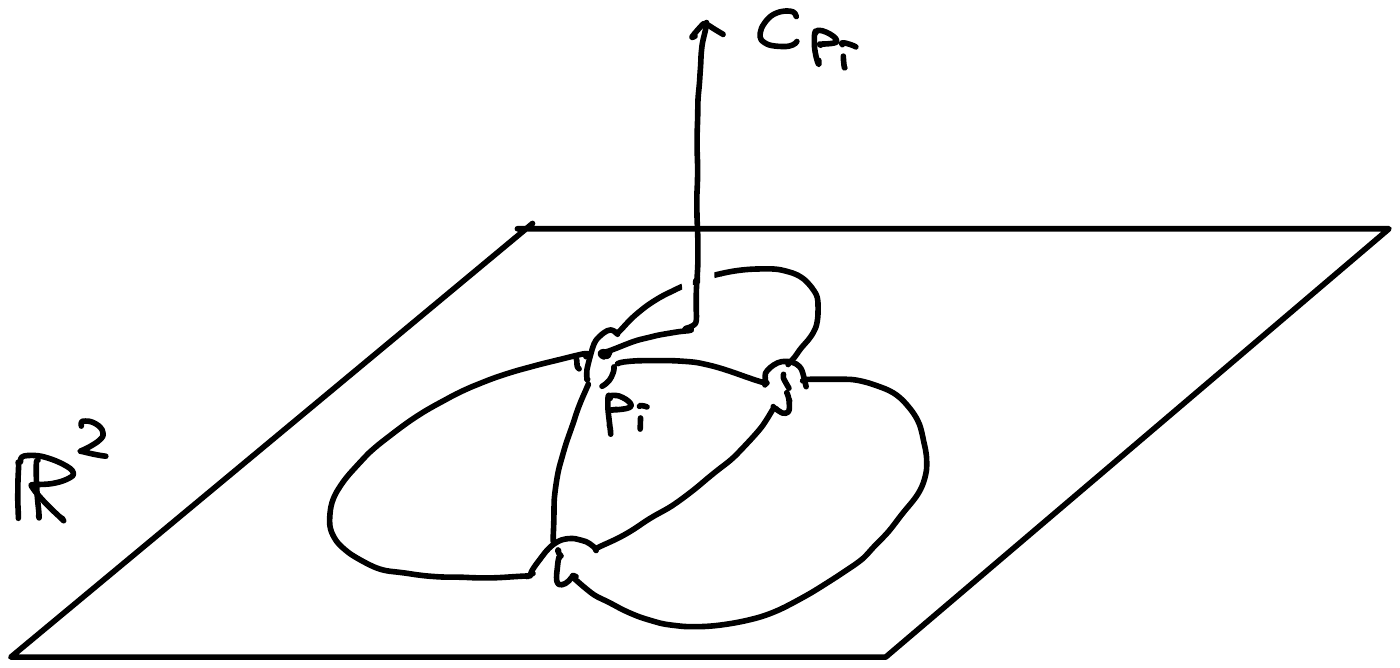}\label{pict912}

For $j=1,\ldots,k+1$, the path $c_{\widetilde{Q}_j}$ parallels to the vector $(0,0,1)$ and is the upward direction. Namely, $c_{\widetilde{Q}_j}=\{\widetilde{Q}_j+s(0,0,1)\mid 0\le s\}$.

\includegraphics[width=8cm]{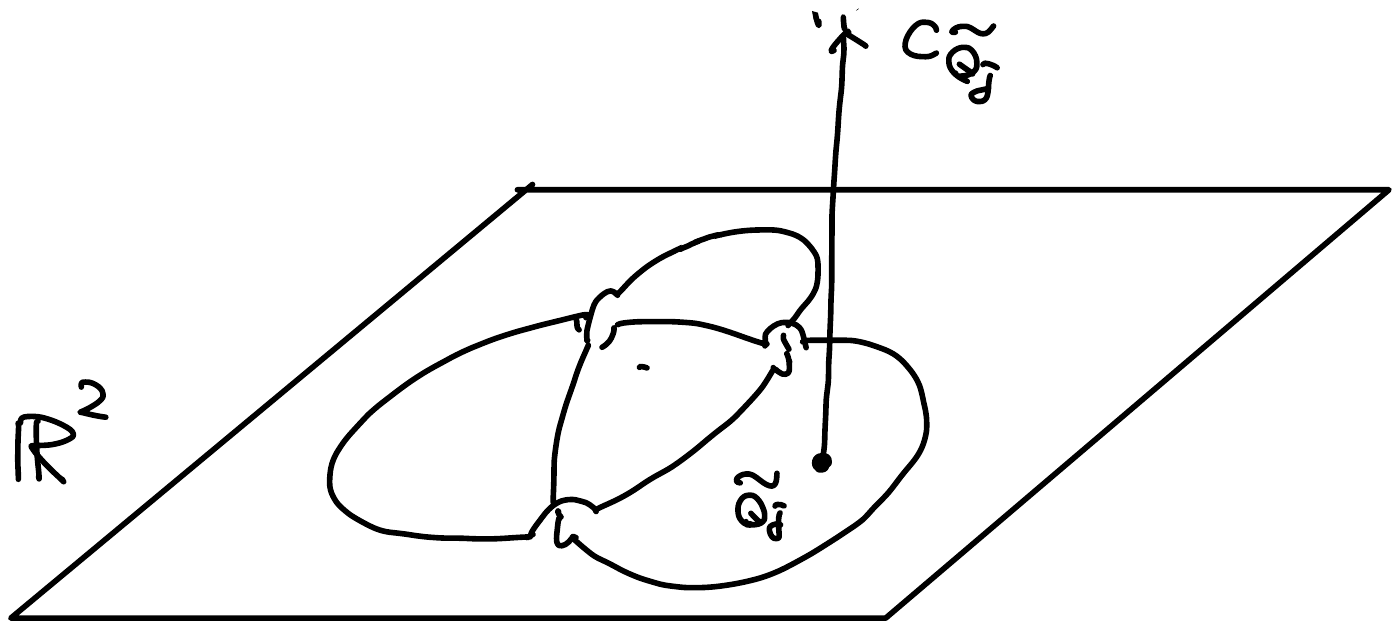}\label{pict913}

The path $c_{\infty}$ connecting $\infty$ and $\infty$ is a trivial path.
\subsection{The Dehn graph $\Gamma_{f,c_f}$ of $f,c_f$ coincides with $\Gamma_K$}
We show that  $\Gamma_K$ coincides with the Dehn graph $\Gamma_{f,c_f}$ of $f$ and $c_f$.
Via the correspondence $p_i\mapsto P_i, q_j\mapsto \widetilde{Q}_j, \infty\mapsto \infty$, we have $V(\Gamma_K)=V(\Gamma_{f,c_f})$.
We have to check that $E(\Gamma_K)(p,p')=E(\Gamma_{f,c_f})(p,p')$ including these labels for any $p,p'\in V(\Gamma_K)=V(\Gamma_{f,c_f})$.
\vskip3mm
\noindent
\underline{Trajectories from $P_i$ to $\widetilde{Q}_j$}
\vskip2mm
There is a one-to-one corresponding between 
$\mathcal M(P_i,\widetilde{Q}_j)$ and $\mathcal D_{P_i}\cap \mathcal A_{\widetilde{Q}_j}\cap f^{-1}(-1.5)$ as sets.
Here $\mathcal D_{P_i}$ is the descending manifold of $P_i$
and $\mathcal A_{\widetilde{Q}_j}$ is the ascending manifold of 
$\widetilde Q_j$ with respect to $-\grad f$.
 
Since $f^{-1}(-1.2)$ is the boundary of a regular neighborhood of $K\cup (\cup_{i=1}^kL_i)$, $\mathcal D_{P_i}\cap f^{-1}(-1.2)$ forms a circle around $P_i$ as in the following left picture.
Then $\mathcal D_{P_i}\cap f^{-1}(-1.5)$ is a simple closed curve around $P_i$ described as in the right picture.

\includegraphics[width=12cm]{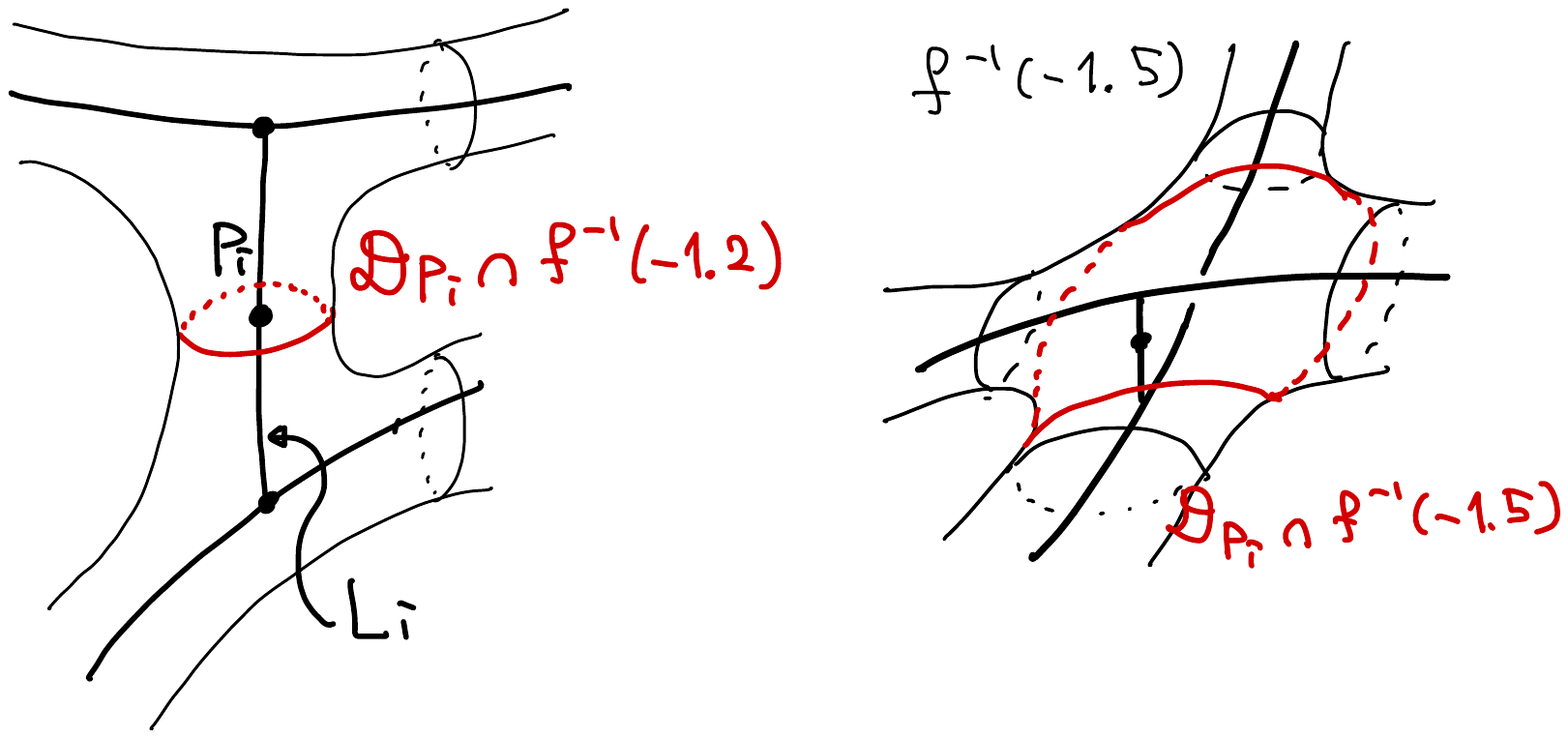}\label{pic_D} 

Similarly,
$\mathcal A_{\widetilde Q_j}\cap f^{-1}(1.5)$ forms a loop around the edge of the region $Q_j$ as in the following picture.

\includegraphics[width=8cm]{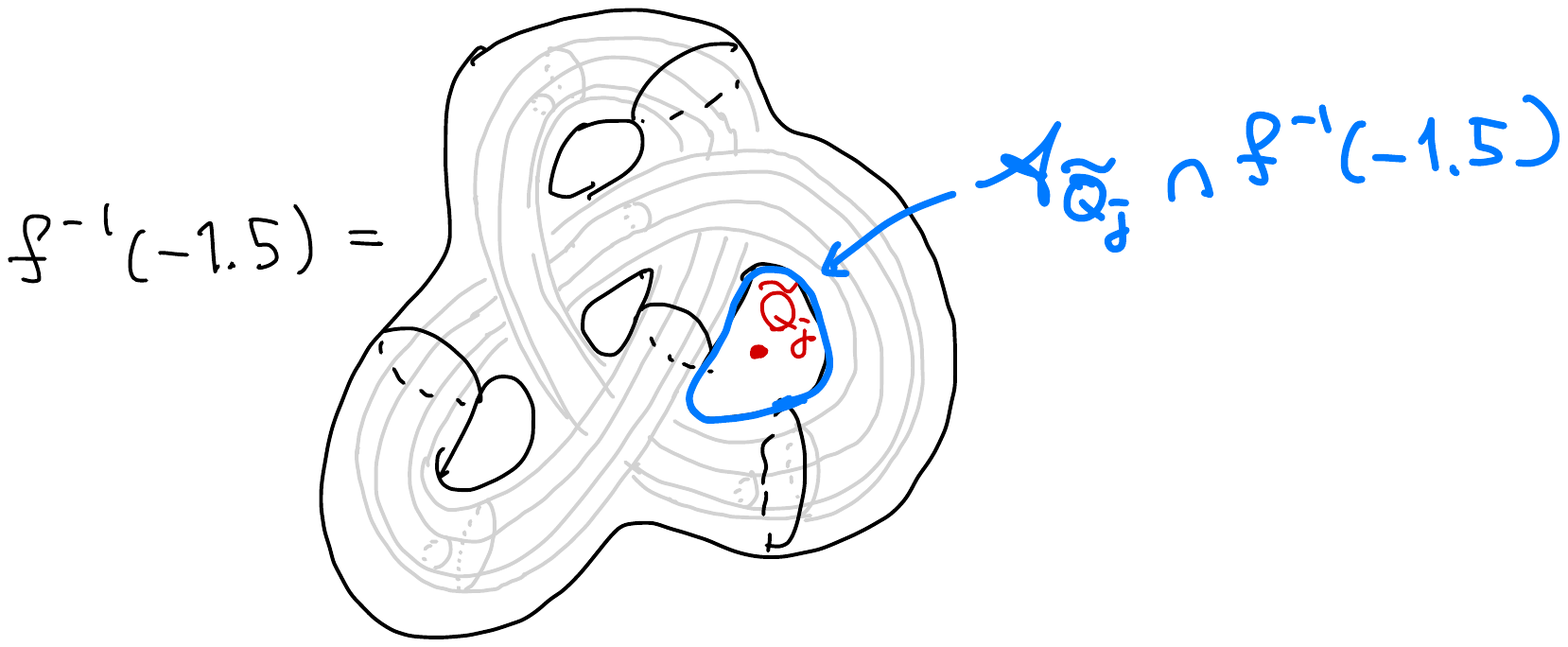}\label{pic_A} 

Therefore we can find out $\mathcal D_{P_i}\cap \mathcal A_{\widetilde Q_j}\cap f^{-1}(-1.5)$ in the corner of the region $P_j$ near $Q_i$.
\vskip3mm
\noindent
{\bf Remark.}
$f^{-1}(-1.5)$ can be considered as a Heegaard surface of a Heegaard splitting compatible with the Morse function $f$.
Actually, a diagram on $f^{-1}(-1.5)$ consisting of several curves occurred from critical points gives a Heegaard diagram of $S^3$.
For example, a similar argument can be used to obtain a doubly pointed Heegaard diagrams that appears in Heegaard Floer homologies. 
See, for example, \cite{OS}, \cite{HeegaardFloer}.
\vskip3mm

Let $\gamma$ be a trajectory from $P_i$ to $\widetilde{Q}_j$.
We next discuss the orientation of $\gamma$ and the class $[c_{\widetilde{Q}_j}\circ \gamma\circ c_{P_i}^{-1}]\in \pi_1(E(K);\infty)$
to compute the label of $\gamma$.

The orientation of each trajectory is determined from that of ascending manifolds, descending manifolds and $E(K)$.
Then the orientation of $\gamma$ depends on the orientation of $\mathcal D_{P_i}\cap f^{-1}(-1.5)$ and that of $\mathcal A_{\widetilde{Q}_j}\cap f^{-1}(-1.5)$.
Recall that these orientations satisfy $T_p{\mathcal A}_p\oplus T_p\mathcal D_p\cong T_pE(K)$ for each critical point $p$. 
We can choose any orientations for $\mathcal D_{P_i}$ and $\mathcal A_{\widetilde{Q}_j}$ by changing the orientations of $\mathcal A_{P_i}$ and $\mathcal D_{\widetilde{Q}_j}$.
Thus we orient $\mathcal A_{\widetilde{Q}_j}\cap f^{-1}(-1.5)$ the counter-clockwise orientation in $\RR^2$ as in the following picture.

\includegraphics[width=8cm]{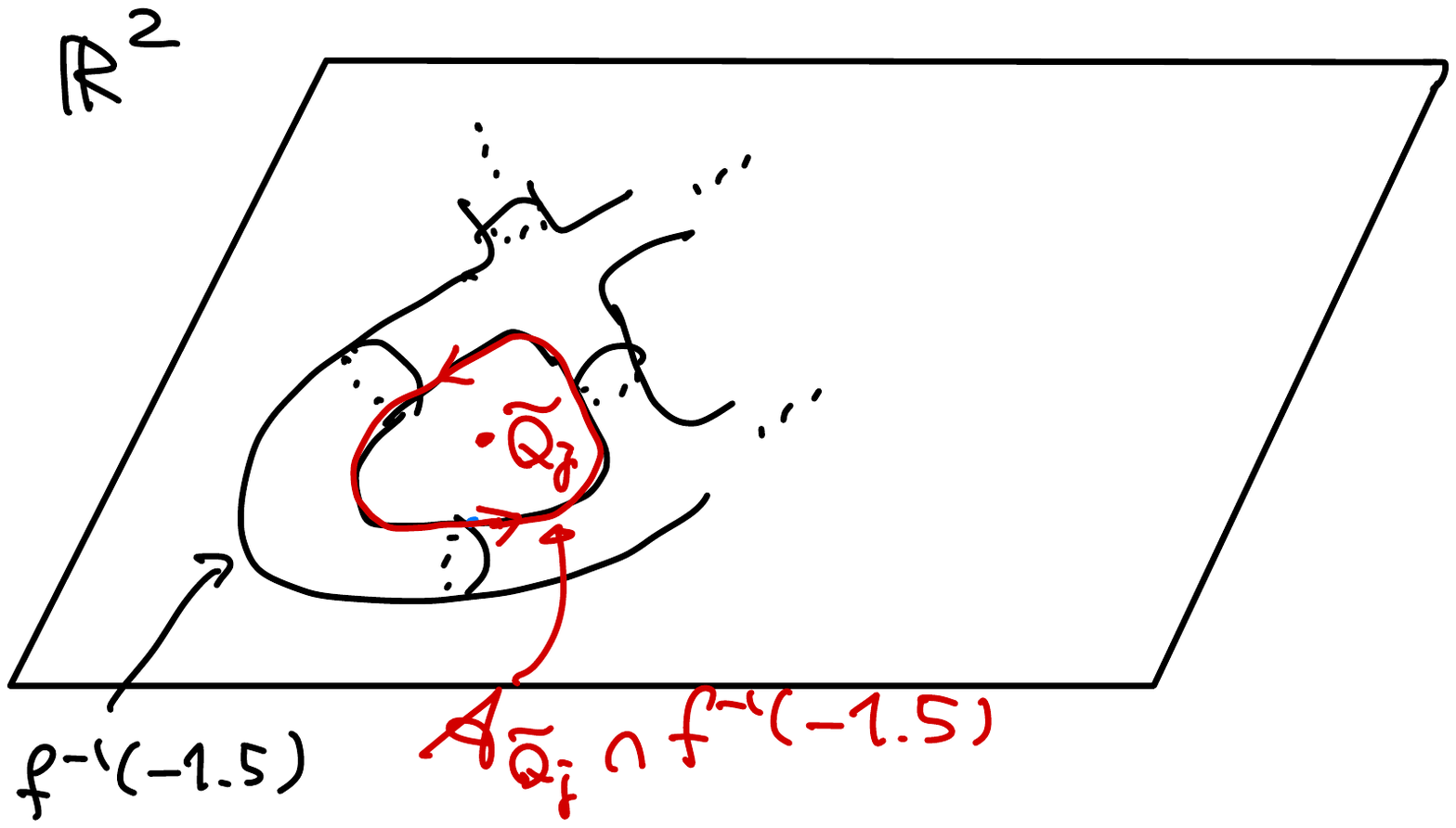}\label{pict_D_ori} 

There are two situation around $\gamma$ as in the following pictures. 
We choose the orientation of $\mathcal D_{P_i}$ so that if $\gamma$ is as in the left picture, $\ep(\gamma)=1$ and if $\gamma$ is as in the right picture, $\ep(\gamma)=-1$.

\begin{center}
\includegraphics[width=10cm]{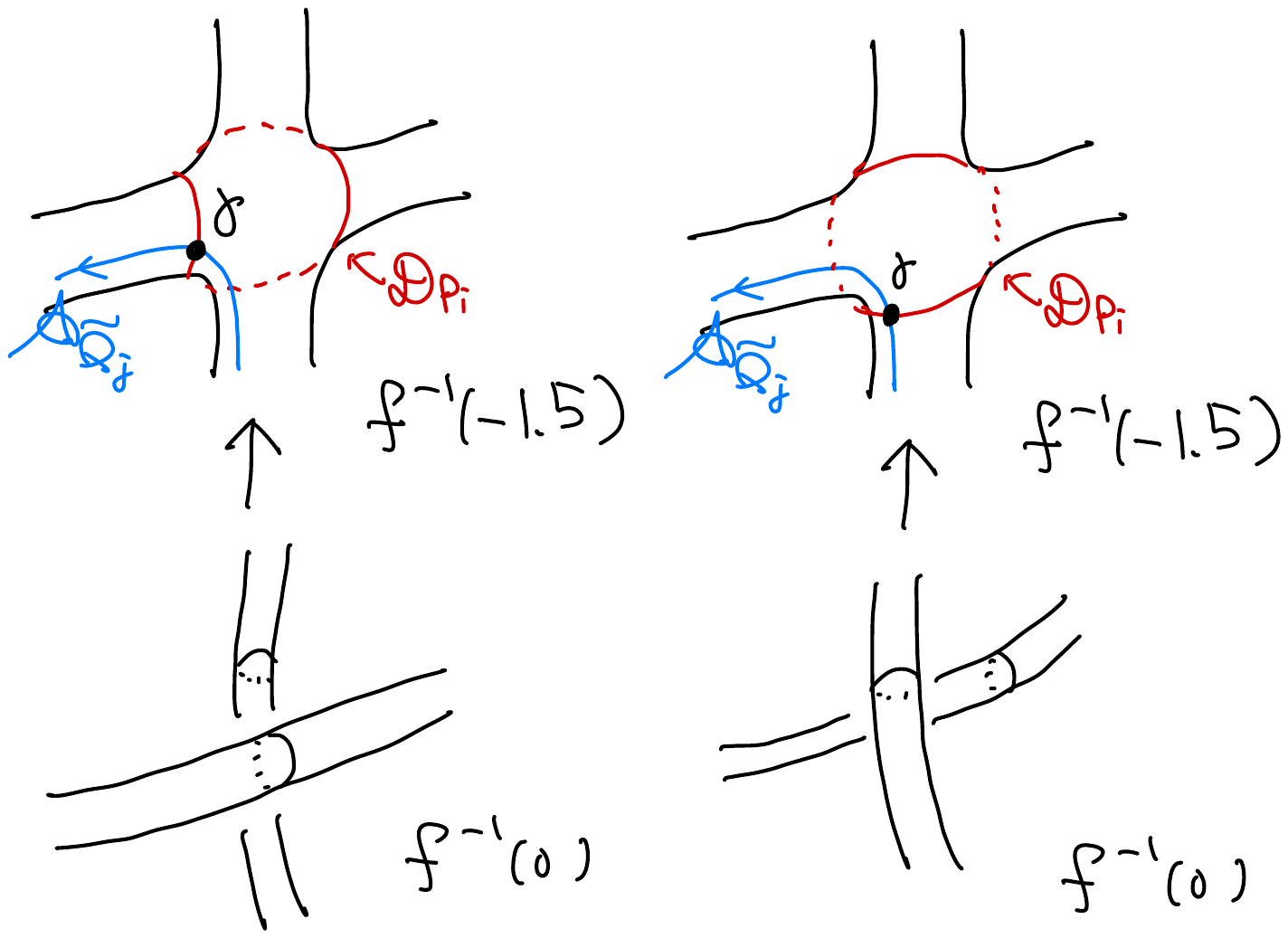}\label{pic_gamma} 
\end{center}

The class $[c_{\widetilde Q_j}\circ\gamma c_{P_i}^{-1}]$ can be  
read from the following picture.

Then the label $l(\gamma)=\ep(\gamma)[c_{\widetilde Q_j}\circ\gamma c_{P_i}^{-1}]$ of $\gamma$ in $\Gamma_{f,c_f}$ equals to the label of $\gamma$ in $\Gamma_K$ given by the rule (R) in Section~\ref{D_1}.

\begin{center}
\includegraphics[width=10cm]{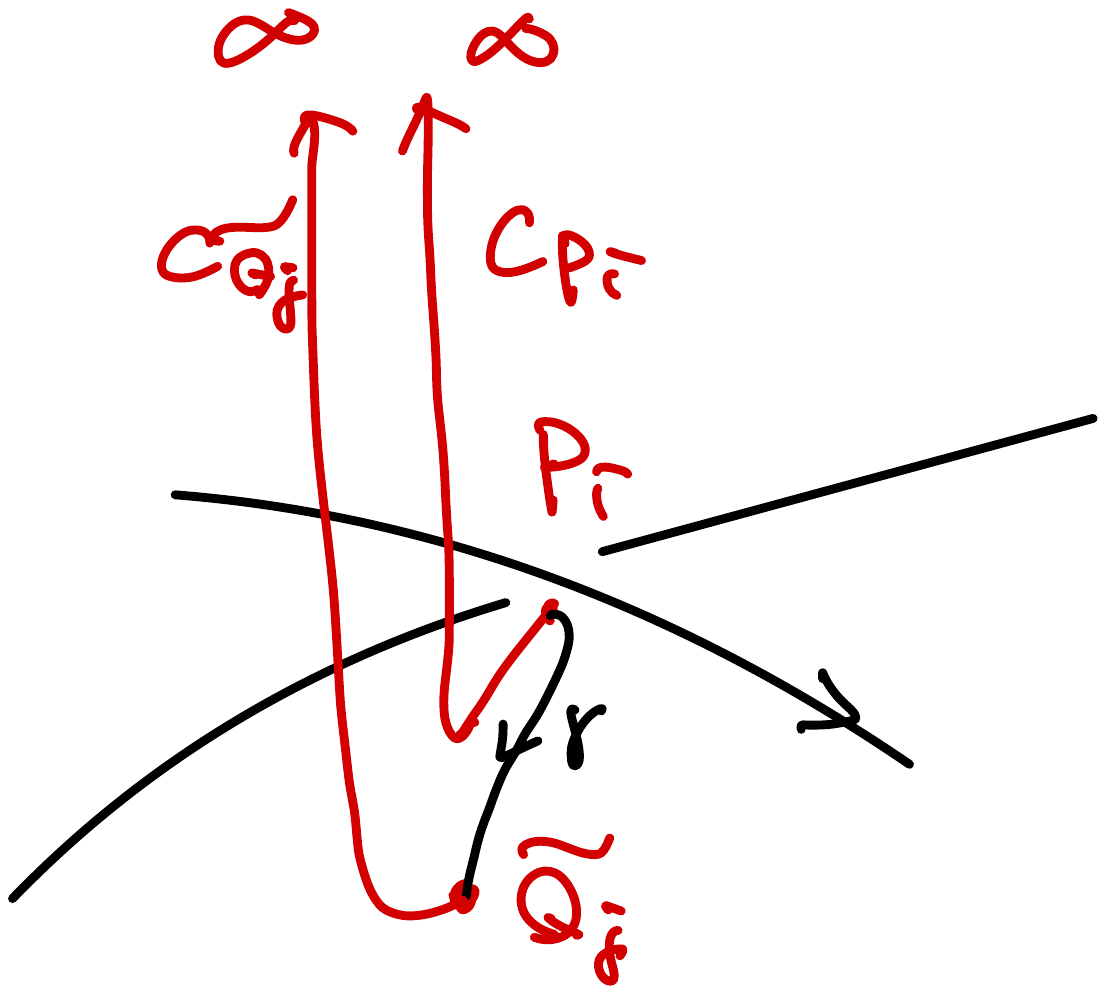}\label{pic_gamma_l} 
\end{center}

\vskip3mm
\noindent
\underline{Trajectories from $\widetilde{Q}_j$ to $\infty$}
\vskip2mm
The descending manifold $\mathcal D_{\widetilde{Q}_j}$ is a 1-manifold parallel to the $z$-axis. 
Then there are two trajectories, denoted by $\gamma_+$ and $\gamma_-$, from $\widetilde{Q}_j$
to $\infty$ for each $j$ as in the following picture.
The orientation of $\mathcal D_{\widetilde{Q}_j}$ given in the above determines orientation of the trajectories as
$\ep(\gamma_+)=1, \ep(\gamma_-)=-1$.
\begin{center}
\includegraphics[width=10cm]{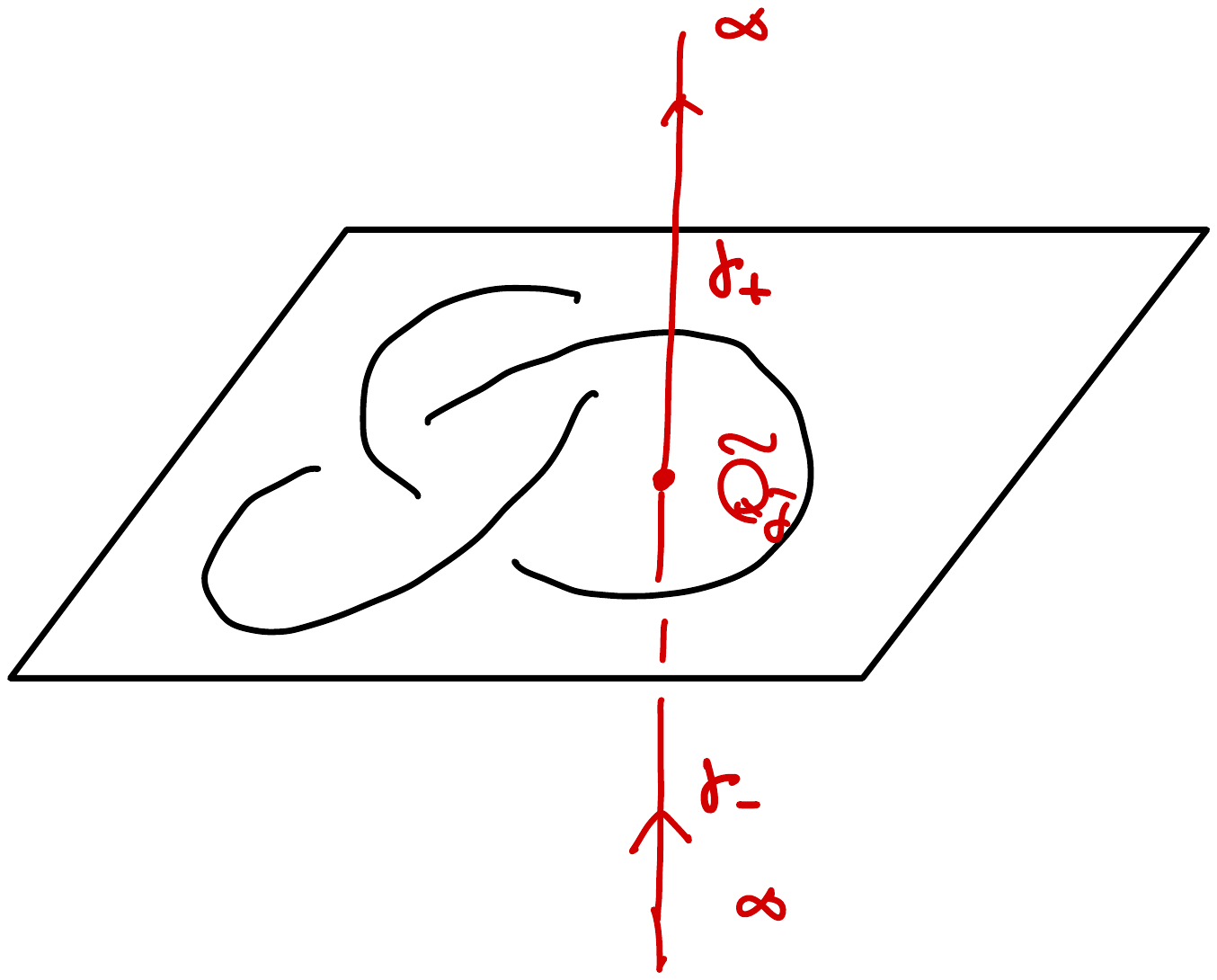}\label{pic_gamma_Q} 
\end{center}
Obviously, $[c_{\infty}\circ \gamma_+\circ c_{\widetilde{Q}_j}^{-1}]
=1\in \pi_1(M,\infty)$. The class 
$[c_{\infty}\circ \gamma_-\circ c_{\widetilde{Q}_j}^{-1}]=
[\gamma_-\circ c_{\widetilde{Q}_j}^{-1}]\in \pi_1(M,\infty)$
coincides with the label assigned to the region $Q_j$ in the diagram $D_2$. 
Then the label $\ep(\gamma_{\pm})[c_{\infty}\circ\gamma_{\pm}\circ c_{\widetilde{Q}_j}^{-1}]$ in $\Gamma_{f,c_f}$ coincides with the label $l(\gamma_{\pm})$ in $\Gamma_K$ given from the diagram $D_2$.

\end{document}